\input amssym.def
\input amssym.tex
\magnification 1200

\def\newline{\break}
\def\item{\vskip .03 in}

\def\min{\setminus}

\def\a{\alpha}
\def\b{\beta}
\def\g{\gamma}

\def\un{\underline}

\def\o{\omega}
\def\O{\Omega}

\def\la{\lambda}

\def\inf{\infty}

\def\longto{\longrightarrow}
\def\mapright#1{\smash
{\mathop{
\longto}\limits^{#1}}}
\def\mapdown#1{\Big\downarrow
\rlap{$\vcenter
{\hbox{$\scriptstyle#1$}}$}}

\def\na{\nabla}
\def\pa{\partial}

\def\bu{\bullet}

\def\O{\underline\Omega}

\def\BR{{ \Bbb R}}

\def\BC{{ \Bbb C}}

\hoffset= 2pc
\hsize=30pc
\vsize=45pc

 \font\sectionfont=cmbx12 at 12pt



\def\sec#1{\vskip 1.5pc\noindent
{\hbox 
{{\sectionfont #1}}}\vskip 1pc}

\def\theo#1#2{\vskip 1pc
  {\bf   Theorem\ #1.} {\it #2}
\vskip 1pc}
\def\lem#1#2{\vskip 1pc{\bf 
\ Lemma #1.}
  {\it #2}\vskip 1pc}
\def\cor#1#2{\vskip 1pc{\bf 
\ Corollary\ #1.}
  {\it #2}\vskip 1pc}
\def\prop#1#2{\vskip 1pc{\bf 
 Proposition \ #1.} 
{\it #2}}

\def\und{\underbar}

\def\uBC*M{\u{\BC}^*~_M}

\def\uBC{\u{\BC}}

\def\calO{{\cal O}}

\def\calC{{\cal C}}

\def\calL{{\cal L}}

\def\text{\textstyle}
\def\pha{\phantom}

\def\text{\textstyle}
\def\pha{\phantom}

\def\2pii{2\pi\sqrt{-1}}

\def\({\left(}
\def\){\right)}
\def\[{\left[}
\def\]{\right]}

\def\NSF92-95{\footnote*{This research was supported in part by NSF grants
DMS-9203517 and DMS-9504522.}}
\def\NSF95{\footnote*{This research was supported in part by NSF grant
 DMS-9504522.}}

\def\328{Lecture Notes in Math.}
\rightline{Revised in February, 1997}
\vskip .124 in
\centerline{\bf {The outer derivation of a complex Poisson manifold.}}
\vskip 1pc

\centerline{by J-L. Brylinski\NSF95 and G. Zuckerman}
\vskip 1.3pc

In this note we introduce the canonical outer Poisson vector field
on a Poisson manifold (both in the smooth and complex-analytic cases).
This means that we construct a canonical class of Poisson vector fields
modulo the locally hamiltonian vector fields. This class is an obstruction
to the existence of a volume form invariant under all hamiltonian flows.
It is defined even both in the case of a smooth Poisson manifold and in the
case of a holomorphic Poisson manifold, even if there is
no global volume form.
We show (in the complex-analytic case) that this outer vector field
 gives a  class of holomorphic functions with singularities
along the complement of the regular set, which is canonical modulo the addition
of Casimir functions. We study this class in terms of various sheaf cohomology
groups of interest.
The same outer vector field has been introduced independently by
Alan Weinstein [{\bf W}]. He has a description of the outer vector field
which relates it to the KMS theory and thus to the Connes-Takesaki theory
of modular automorphisms of a von Neumann algebra.

We thank Alan Weinstein for communicating to us an early version
of his work and for useful correspondence. The first author thanks
the participants at the Singularity Seminar at Marseille-Ch\^ateau
Gombert where he lectured on this work. We thank Philip Foth for useful
comments on a previous version of this note.

\sec{1. Hamiltonian  and outer vector fields 
of a Poisson manifold}

We will work either with a $C^{\inf}$ or a complex
Poisson manifold $M$. In the first case there is a Poisson tensor
$\pi\in \wedge^2TM$, where $TM$ is the tangent bundle.
In the second case the Poisson tensor $\pi$ is a holomorphic section
of $\wedge^2\Theta_M$, where $\Theta_M$ is the holomorphic
tangent bundle. 

There are three interesting classes of vector fields, which
we enumerate starting with the largest class:

(1) the {\it Poisson vector fields}: a vector field $\xi$
is Poisson it it preserves the Poisson structure, that
is $$\calL_{\xi}\pi=\{\xi,\pi\}=0,$$ where $\{\pha{a},\pha{a}\}$
is the Schouten bracket.

(2) the {\it locally hamiltonian vector fields}: locally
$\xi$ is of the form $X_H=i(dH)\pi$, where $i$ denotes interior
product, and $H$ is a smooth (resp. holomorphic)
function defined locally.

(3) the hamiltonian vector fields: $\xi=X_H$ for some
global $H$.

It is interesting to devise cohomological criteria
separating these types. For instance, let $\xi$ be
 locally hamiltonian; then there is an open covering $(U_i)$
such that $\xi=X_{H_i}$ over $U_i$. Then $f_{ij}=H_j-H_i$
satisfies $X_{f_{ij}}=0$, so it is a Casimir function over $U_{ij}$.
We have a \u Cech $1$-cocycle $f_{ij}$ with values in the sheaf
$\calC_M$ of smooth (resp. holomorphic) Casimir functions.

\prop{1}{(1) For $M$ a smooth Poisson manifold, we have 
an exact sequence
$$0\to  Ham(M)\to Loc~ Ham (M)\to H^1(M,\calC_M)\to 0.$$
\item
(2) For $M$ a complex Poisson manifold, we have 
an exact sequence
$$0\to  Ham(M)\to Loc~ Ham (M)\to H^1(M,\calC_M)\to H^1(M,\calO_M),$$
where $\calO_M$ is the sheaf of germs of holomorphic functions.}

Now we wish to examine the difference between the first two types
of vector fields. We need to sheafify the question. Let
$\un{ Ham}$ be the sheaf of locally hamiltonian vector fields, and
let $\un{ Der}$ be the sheaf of Poisson vector fields. Then
$\un{ Ham}$ is a subsheaf of $\un {Der}$. If we have a vector field
$\xi$ which is Poisson, it gives a section of the quotient
sheaf $\un{ Der}/\un{ Ham}$. This quotient sheaf is hard to describe in
general. In the regular case we have the following description, where
$\un{TX}$denotes the sheaf of smooth tangent vectors, 
$\un{\Theta_X}$ the sheaf of holomorphic tangent vectors.

\prop{2}{Let $M$ be a regular Poisson manifold
such that there exists a smooth fibration $p:M\to X$ 
whose fibers are the symplectic leaves.
Then the quotient sheaf $\un{ Der}/\un{ Ham}$ identifies with
the pull-back sheaf $p^{-1}\un{ TX}$ (in the smooth
case), resp. $p^{-1}\un{ \Theta_X}$ (in the holomorphic case).}

Now in the smooth case the sheaf
 cohomology $H^i(M,p^{-1}\un {TX})$
is equal to $Vect(X)\otimes H^i(F,\BR)$, where
$F$ is the fiber of $p$. In the holomorphic case, if $M$ is Stein, the group 
$H^i(M,p^{-1}\un{\Theta_X})$ is equal to $\Theta(X)\otimes H^i(F,\BC)$.
In general there is a spectral sequence with $E_2$ term
$E_2^{pq}=H^p(X,\un{\Theta_X}\otimes H^q(F,\BC))$
converging to $H^i(M,p^{-1}\un{\Theta_X})$.
So we obtain

\prop{3} {Let $M$ be as in Proposition 2.
Let $\xi$ be a Poisson vector field on $M$. Then
$\xi$ is locally hamiltonian iff
the corresponding vector field on $X$ is trivial.}

Now we come to the notion of an outer Poisson vector field.
This is defined to be a section of the quotient sheaf
$\un{ Der}/\un{ Ham}$.

\lem{1}{Let $M$ be a smooth Poisson manifold or a complex Poisson
manifold which is Stein. An outer Poisson vector field comes from a global
Poisson vector field iff the corresponding class in
$H^1(M,\un{Ham})\tilde{\to}H^2(M,\calC_M)$ is trivial.}

The isomorphism comes from the exact sequences of sheaves
$$0\to\calC_M\to \un{C}^{\inf}_M\to \un{Ham}_M\to 0,$$
in the smooth case
$$0\to\calC_M\to \calO_M\to \un{Ham}_M\to 0,$$
in the holomorphic case. In the holomorphic case we
have $H^1(M,\calO_M)=H^2(M,\calO_M)=0$ because $M$ is Stein.

\sec{2. Koszul operators and the canonical outer vector field.}

The results here were also obtained independently by A. Weinstein
[{\bf W}]. Following Koszul [{\bf K}], we say that an operator
$D:\wedge^iTM\to\wedge^{i-1}TM$ generates the Schouten bracket iff
we have
$$\{ u,v\}=(-1)^p[D(u\wedge v)-(Du)\wedge v-(-1)^pu\wedge(Dv)],\eqno(1)$$
where $u\in\wedge^pTM$, $v\in\wedge^qTM$. There is a similar
notion in the holomorphic case. We call these operators
\und{ Koszul operators}.

\theo{1} {(Koszul) (1) The set of Koszul operators $D$ is in canonical
bijection with the set of connections $\na$ on the canonical bundle
$\wedge^nT^*M$ (resp. $\o_M=\O^n_M$).
\item
(2) The square $D^2$ is interior product with the curvature
$K$ of $\na$.
\item
(3) Let $\na$ be a connection on the canonical bundle, and
let $\a$ be a $1$-form. Then the Koszul operator
corresponding to the connection $\na+\a$ is $D+i(\a)$.}

In particular, a Koszul operator is of square $0$ iff
the corresponding connection is flat.

Recall that a holomorphic vector bundle $E\to M$ has Chern
classes $c_p(E)\in H^p(M,\O^p)$ in Hodge
cohomology, and slightly more refined classes
in $H^p(M,\O^p_{M,cl})$, where $\O^p_{M,cl}$ is the sheaf
of closed holomorphic $p$-forms.

\cor{}{Let $M$ be a complex manifold. 
\item
(1) There exists a global holomorphic
Koszul operator iff the first Chern class
$c_1(\Theta_M)\in H^1(M,\O^1_M)$ of $\Theta_M$ 
in Hodge cohomology is $0$.
\item
(2)
There exists a global holomorphic Koszul operator of square $0$ iff the 
class $c_1(\Theta_M)\in H^1(M,\O^1_{M,cl})$ is $0$.}

We recall how a connection $\na$ defines a Koszul operator.
We discuss only the holomorphic case.
We have the canonical isomorphism
$$*:\wedge^i\Theta_M\tilde{\to}\O^{n-i}_M\otimes\o_M^{\otimes-1},$$
where $\O^i_M$ is the sheaf of holomorphic $i$-forms.
The connection $\na$ on $\o_M$ induces
one on the dual line bundle $\o_M^{\otimes-1}$, which
in turn extends to operators
$\O^j_M\otimes\o_M^{\otimes-1}\to \O^{j+1}_M\otimes\o_M^{\otimes-1}$.
Using $*$ this induces an operator
$D:\wedge^i\Theta_M\to\wedge^{i-1}\Theta_M$. It is easy to see
that $D$ is a Koszul operator.

Note that any volume form $\nu$ induces a flat connection $\na$ on $\o_M$, hence
a Koszul operator of square $0$ which we will denote by $D_{\na}$. More generally any
multi-valued volume form $\nu$, such that monodromy around any loop multiplies
$\nu$ by some constant, also defines a flat connection.

Now assume that $M$ is  a Poisson manifold with Poisson tensor $\pi$.
Then for any (local) Koszul operator $D$ of square $0$, the vector field $D\pi$ is Poisson
because (1) gives
$$0=D\{\pi,\pi\}=-(D\pi)\wedge\pi-\pi\wedge D\pi=-2(D\pi)\wedge\pi$$
and then
$$
2\{ D\pi,\pi\}=2D^2(\pi)\wedge\pi-D^2(\wedge^2\pi)=0$$
(cf. [{\bf K}]).

Now if we change the Koszul operator of square $0$   to $D+i(df)$,
then we change $D\pi$ to $D\pi+X_f$. Hence
we have

\theo{2}{There exists a canonical section $\xi$ of $\un{Der}/\un{Ham}$,
which is locally equal to $D\pi$ for any Koszul operator $D$ of square $0$.}

This is called the {\it canonical outer vector field}.

The  canonical outer vector field can be computed more concretely from any local
volume form $\nu$. We have

\prop{4} {Given a volume form $\nu$ over $U\subseteq M$, a representative
of $\xi$ over $U$ is given by the  vector field $v$ such that
$(v\cdot f)\nu=\calL_{X_f}\nu$.}

Let us now look at the regular case. Then we have

\prop{5} {For a regular Poisson manifold equipped with a fibration
$p:M\to Z$ whose fibers are the symplectic leaves, and for any volume
form $\nu$, the  Poisson vector field $D\pi$ is equal
to 
$$\pm X_{log~f},$$
where the function $f$ is defined as follows.
Let $\a$ be some volume form on $Z$; then the Liouville fiberwise volume form
$\b$ defines a volume form $\a\wedge\b$ on $M$. Then
$$f={\a\wedge\b\over \nu}.$$}

\cor{}{(see also [{\bf W}]) The section of $\un{Der}/\un{Ham}$
is supported in the singular set of the Poisson manifold $M$.}

As an example, let $\frak g$ be a complex Lie algebra,
and let $\la:\frak g\to \BC$ be the modular character
$\la(\xi)=Tr(ad(\xi)$.
Then the constant vector field $\la$ on $\frak g^*$
represents the outer vector field $\xi$. We have

\prop{6} {Let $\frak g$ be a (real or complex) Lie
algebra which is not unimodular. Then the outer Poisson vector field
$\xi$ does not vanish in any neighborhood of $0$.}

This is proved easily by looking at the Taylor series at the origin
of  any hamiltonian vector field.

\sec{3. The canonical singularity class.}

For a Poisson manifold $M$, we have the canonical Poisson
outer vector field
$\xi$. We have proved that its restriction to the regular open set $U$
is trivial. We wish to think of $\xi$ as corresponding to a class
(modulo Casimir functions) of singularities along $S=M\min U$
of multivalued meromorphic functions. This will be mostly interesting in the complex-analytic
case. We use cohomology $H^i_S(M,F)$  with supports in $S$ of a sheaf $F$, 
which fits into the exact sequence
$$0\to H^0_S(M,F)\to H^0(M,F)\to H^0(U,F)\mapright{\pa}
H^1_S(M,F)\to H^1(M,F)\to\cdots.$$

\prop{7}{(1) Assume that $M$ is equipped with a holomorphic volume
form $\nu$, and let $D$ be the corresponding Koszul operator.
Then the Poisson vector field $v=D\pi$ is locally hamiltonian
over $U$, hence it gives rise to a class $\bar f_{\nu}$
in $H^0(U,\calO_M/\calC_M)$.
\item
(2) The outer vector
 field $\xi$ gives a canonical class
in $H^1_S(M,\calO_M/\calC_M)$.
If there is a global volume form, this class
is the image of the class $\bar f_{\nu}$ in (1) by the boundary operator $\pa$.}

Here $H^1_S(M,\calO_M/\calC_M)$ denote sheaf cohomology
\und{with supports in $S$.}

In (1) we could more generally assume that $M$ has a multi-valued
volume form $\nu$ such that $T\nu\over\nu$ is constant for any monodromy
transformation $T$.

This class can be understood as follows. First assume there
is a global volume form. The global
 Poisson vector field $v=D\pi$ which over $U$ is locally hamiltonian, 
hence given by a section of $\calO_M/\calC_M$ over $U$.
We can change the volume form, multiplying it by  a non-vanishing function
$f$ over $M$. Then $log~f$ is  globally defined over $U$ modulo constants,
hence it is globally defined as a section
of the quotient sheaf $\calO_M/\calC_M$ over $U$. Therefore the image
in $H^1_S(M,\calO_M/\calC_M)$ is well-defined.

Now relax the assumption about a global volume form.
The cohomology with supports
$H^*_S(M,\calO_M/\calC_M)$ can be computed as a direct limit
over open coverings $(U_i)$ of $M$
of the quotient complex $C^{\bu}((U_i\min (S\cap U_i),F)/C^{\bu}((U_i),F)$
where $F=\calO_M/\calC_M$.
This is a complex whose $(p+1)$-st term is
$$\oplus_{i_0,\cdots,i_p}\Gamma(U_{i_0\cdots i_p}
\min (S\cap U_{i_0\cdots i_p}),F)
/\Gamma(U_{i_0\cdots i_p},F).$$
This description is correct because there is no section
of $F$ on an open set which vanishes on a dense open subset.

Now cover $S$ by open sets $S\cap U_i$, $U_i$ open in $M$,
 such that there is a volume form
$\nu_i$ over $U_i$. 
Then we have a corresponding Poisson vector field
$v_i$ over $U_i$; the restriction of $v_i$ to
$U_i\min(S\cap U_i)$ is locally hamiltonian, hence defines a section
$g_i$ of $\un{Ham}\tilde{\to}\calO_M/\calC_M$ over this open set.
Now we have: $v_j-v_i=X_{log{\nu_j\over\nu_i}}$ over $U_{ij}$
hence if we put $f_{ij}=log~{\nu_j\over\nu_i}$, we get a section
of $\calO_M/\calC_M$ over $U_{ij}$ such that
$v_j-v_i=X_{f_{ij}}$, hence $g_j-g_i=f_{ij}$. It follows
that the \v Cech coboundary of $(g_i)$ is zero in the quotient
complex, hence we have indeed a cohomology class in degree $1$.

There is another interesting cohomological invariant attached to $\xi$.
This is the class in $H^2_S(M,\calC_M)$ obtained
from the class in $H^1_S(M,\calO_M/\calC_M)$ by applying the boundary map for the exact
sequence of sheaves
$0\to \calC_M\to \calO_M\to \calO_M/\calC_M\to 0.$

If we have a volume form on $M$, then there is a corresponding class
in \hfill\break$H^1(U,\calC_M)$. We can summarize the situation into the following square
$$\matrix{H^0(U,\calO_M/\calC_M)&\mapright{\pa}&H^1(U,\calC_M)\cr
\mapdown{\pa}&&\mapdown{\pa}\cr
H^1_S(M,\calO_M/\calC_M)&\mapright{\pa}&H^2_S(M,\calC_M).}$$
The classes in the groups of the first row depend on a volume form;
the classes in the groups in the second row are entirely intrinsic.

We would like to suggest that the groups $H^1(U,\calC_M)$ and
$H^2_S(M,\calC_M)$ have a geometric interpretation in terms
of multivalued holomorphic functions on $U$ whose variation
(under monodromy) is a Casimir function. We first study the case
where the Poisson tensor $\pi$ is non-degenerate on a dense open set
$U$. Then the vanishing locus of the section $\wedge^n\pi$
of $\wedge^n\Theta_X$ is a divisor $D=\sum_i~n_iD_i$. The singular locus
$S$ is the support of $D$, i.e., the union of the $D_i$.
Since all the Casimir functions on any open set are constant,
the sheaf $\calC_M$ is the constant sheaf $\un\BC_M$.
We have $H^2_S(M,\BC)=\sum_i~\BC[D_i]$, where $D_i$ is the class
of the divisor  in cohomology with support. Then we can state

\prop{8}{The class in $H^2_S(M,\BC)$ corresponding to
 the canonical outer vector field
is the class $\sum_i~n_i[D_i]$ of the divisor of $\wedge^n\pi$.}

We illustrate this in a simple
example. Let $\frak g$ be the two-dimensional
complex Lie algebra
with basis $x,y$, such that $[x,y]=x$. Then $\frak g^*=\BC^2$ is a holomorphic
Poisson manifold, with Poisson tensor $\pi=x{\pa\over \pa x}\wedge{\pa\over\pa y}$.
The regular open set $U$ is the open set $U=\{(x,y);x\neq 0\}$.
The singular set $S$ is the line $x=0$. If we use the volume form $dx\wedge dy$,
then we find that $\xi$ is represented by the Poisson vector field
$v={\pa\over \pa y}$. Over a contractible open set in $U$, this is the hamiltonian
vector field associated to a branch $\log~x$ of the logarithm function.
The obstruction to finding a global function $f$ such that $v=X_f$ 
is given by the monodromy
of $\log~x$, which is encoded in the generator of $H^1(U,\BC)=\BC$.
The class in $H^2_S(M,\BC)=\BC$ is the canonical orientation class.

There is a similar computation for quadratic Poisson algebras, and for the
Poisson-Lie group corresponding to $SL(2)$.

In general the group $H^1(U,\calC_M)$ is hard to describe. However
for any symplectic leaf $F\subset U$ we have a restriction map
$H^1(U,\calC_M)\to H^1(F,\BC)$; therefore for any loop $\g\in F$ we have
a corresponding character $\chi_{\g}:H^1(U,\calC_M)\to\BC$. This has
the following description:

\prop{9}{Let $U$ be an open subset of $M$ over which there is a fibration
$p:M\to Z$ whose fibers are symplectic leaves. Assume there is a 
multi-valued global volume form $\nu$ on $U$, and let
$c\in H^1(U,\calC_M)$ be the corresponding cohomology class. Let
$F\subset U$ be a symplectic leaf, and let Let $x=p(F)$. Let
$\a$ be an element of $\wedge^{max}\Theta_{X,x}^*$, and
let $g$ be the multi-valued function on $F$ given by  
$$g={\a\wedge\b\over\nu},$$
where $\b$ is the fiberwise volume form 
as in Proposition 5. Then for any loop $\g\in F$ the
number $\chi(c)$ is equal to the variation of the multi-valued
function $\log~g$ along the loop $\g$.}

Concerning the cohomology groups $H^i(M,\calC_M)$ associated to 
a complex Poisson
manifold $M$, they are modules
over the algebra $H^0(M,\calC_M)$ of global Casimir functions.
It is tempting to conjecture that for $M$ a Stein manifold, these are 
finitely-presented modules.
\vskip 1pc

{\bf References}:
\vskip 1pc

[{\bf K}] J-L. Koszul, \it Crochet de Schouten-Nijenhuis et cohomologie, \rm
Ast\'erisque vol. \bf Hors-S\'erie\rm: Elie Cartan et les Math\'ematiques
d'Aujourd'hui, Soc. Math. Fr. (1985), 257-271

[{\bf W}] A. Weinstein, \it The modular automorphism group of
a Poisson manifold, \rm preprint (1995)

\bye